\documentstyle[12pt, leqno]{amsart}
\input{amssym.def}
\input{amssym.tex}

\newtheorem{theorem}{Theorem}[section]

\theoremstyle{definition}

\newtheorem{remark}[theorem]{Remark}
\newtheorem{claim}[theorem]{Claim}

\def\Cal{\cal}

\def\per{{\rm{Per}}}

\def\cS{{\Cal S}}

\font\bb=msbm10 at 12pt
\font\bsub=msbm10 at 10pt
\font\bsubsub=msbm10 at 6pt
\newcommand{\bbR}{\mbox{\bb R}}

\newcommand{\bbZ}{\mbox{\bb Z}}
\newcommand{\bbQ}{\mbox{\bb Q}}
\newcommand{\bbN}{\mbox{\bb N}}
\newcommand{\bbT}{\mbox{\bb T}}
\newcommand{\bbA}{\mbox{\bb A}}

\newcommand{\bsubZ}{\mbox{\bsub Z}}

\newcommand{\bsubsubR}{\mbox{\bsubsub R}}
\begin{document}

\title[Dynamical systems arising from elliptic curves]{Dynamical systems
arising\\ from elliptic curves}
\author[D'Ambros \and Everest \and Miles \and Ward]
{P. D'Ambros \and G. Everest \and R. Miles \and T. Ward}
\address{School of Mathematics, University of East Anglia,
Norwich NR4 7TJ, UK}
\email{g.everest@@uea.ac.uk}
\dedicatory{Dedicated to the memory of Professor Anzelm Iwanik}
\thanks{The first author gratefully acknowledges
the support of INdAM, the
third of E.P.S.R.C. grant 97700813}

\begin{abstract}
We exhibit a family of dynamical systems arising from rational
points on elliptic curves in an attempt to mimic the
familiar toral automorphisms. At the non-archimedean primes, a
continuous map is constructed on the local elliptic curve whose topological
entropy is given by the local canonical height. Also, a precise 
formula for the periodic points is given. There follows a 
discussion of how these local results may be glued together to give
a map on the adelic curve.  We are able to give a map whose entropy is the global canonical height
and whose periodic points are counted asymptotically by the
real division polynomial (although the archimedean component of the
map is artificial). Finally, we set out a precise conjecture
about the existence of elliptic dynamical systems and discuss a possible
connection with mathematical physics.
\end{abstract}

\maketitle

\section{Introduction}

Let $F\in \bbZ [x]$ denote a primitive polynomial with degree $d$, which 
factorizes as
$F(x)=b\prod_i(x-\alpha_i)$. Then $F$ induces a homeomorphism
$T_F$ on a compact, $d$-dimensional group $X=X_F$, via the companion
matrix of $F$. The group $X$ is an example
of a {\sl solenoid} whose definition is discussed in Section 2 below. The essential properties
of this dynamical system $T_F:X\rightarrow X$ are as follows.
\begin{enumerate}
\item The topological
entropy $h(T_F)$ is equal to $m(F)$, the Mahler measure of $F$ (see
(\ref{mahlermeasure}) below).
\item Let $\per_{n}(F)$ denote the subgroup of $X$ consisting of elements
of period $n$ under $T_F$, $\per_{n}(T_F)=\{x\in X:T_F^n(x)=x\}$.
If no $\alpha_i$ is a root of unity, then
$\per_{n}(T_F)$ is finite with order
$
|\per_{n}(T_F)|=|b|^n\prod_i|\alpha_i^n-1|.
$
\end{enumerate}

For background, and proofs of these statements, see \cite[II]{everest-ward-1999}
and \cite{schmidt-1995-algebraic-dynamical-systems}.
The solenoid is a generalization of the (additive) 
circle denoted $\bbT$. Indeed, 
when $F$ is monic with constant
coefficient $\pm1$,
we may take $X$ to be $\bbT^d$; the
resulting map is the automorphism of the
torus determined by the companion matrix to
$F$. The immanence of the 
circle
is seen in both 1 and 2 above. In 1, we may take the definition
of the Mahler measure to be the logarithmic integral of $|F|$ over the
circle,
\begin{equation}\label{mahlermeasure}
m(F)=\int_0^1 \log |F(e^{2\pi it})|dt=\log |b|+\sum_{i=1}^d\log ^
+|\alpha_i|,
\end{equation}
where the last equality follows from Jensen's Formula
(see \cite[Lemma 1.9]{everest-ward-1999}
for proof).
In 2, the periodic points formula is equivalent to evaluating the
$n$-th division polynomial of the circle on the zeros of $F$. That is,
if we take $\phi_n(x)=\prod_{\zeta^n=1}(x-\zeta)$, we get the 
formula
$\alpha^n-1=\phi_n(\alpha)$,
so
$\vert\per_n(T_F)\vert=\vert b^n\times\prod_{i}\phi_n(\alpha_i)\vert$.

In \cite{everest-ward-1998} and \cite[IV]{everest-ward-1999} we set
out reasons for believing in the existence of elliptic dynamical
systems, regarding properties 1 and 2 as the paradigm. Assuming $d=1$ for
example, there ought to be a dynamical system where the immanent group
is a rational elliptic curve with the zero of $F$ corresponding
to the $x$-coordinate of a rational point on that curve. In other words, for every elliptic curve $E$ and every point $Q \in E(\bbQ)$
we are seeking a continuous map $T=T_Q:X\rightarrow
X$ on some compact space $X=X(E)$ whose dynamical data should be described by
well known quantities associated to the point on the curve. We expect the
entropy of $T$ to be the global 
canonical height of the point $Q$ (a well-known analogue of Mahler's
measure) and the
elements of period $n$ should be related to the elliptic
$n$-th division polynomial evaluated 
at the point $Q$. 

There now follows a brief description of this paper, explaining where to look
for our main conclusions. For reasons we will present in Sections 2 and 3, 
it is to be expected
that the underlying space $X$ should be the adelic curve. Section 2 
recalls the classical definition of
the solenoid and the action $F$ induces
on it. Lind and Ward
\cite{lind-ward-1988} re-worked the classical 
theory in adelic
terms. They showed that the topological
entropy can be decomposed into a sum
of local factors, each of which is the entropy of a corresponding
local action. Each of these local factors can be identified as a
corresponding local component of the Mahler measure. Section 3 recalls the
basic theory of elliptic curves needed. In particular, the decomposition of the global canonical height into
a sum of local factors. Also, we recall that the $p$-adic curve
is isomorphic to a simpler group, on which we may expect
to define dynamical systems. In Section 4, in particular the
conclusion,
we will construct a dynamical system where the underlying space
is a $p$-adic elliptic curve and where the map is induced by a point
on that curve. The map in question is a 
$p$-adic analogue of
the well-known $\beta$-transformation. The entropy of the map is the
local canonical height of the point, and the periodic points can
be counted exactly. In Section 5, we will 
consider how to glue together these local maps to get a
global dynamical system. Here we are given an elliptic curve $E$
defined over $\bbQ$ and a rational point $Q \in E(\bbQ)$. The point $Q$
induces a dynamical system where the underlying 
space is the elliptic adeles and the entropy turns out to be the
global canonical height of the point $Q$.
The construction of the map at the archimedean
prime is artificial since it relies upon {\it a priori} knowledge of the
height of the point (although it is a curious coincidence that 
the map is a classical 
$\beta$-transformation). We hope this will bring into better focus the
construction at the non-archimedean
primes where the map uses no such {\it a priori} knowledge of the height.
The artificiality of the map is somewhat redeemed when we go on to show that
the periodic points are counted asymptotically by the real division polynomial
at the point $Q$. This last result makes use of some non-trivial results; 
one from
elliptic transcendence theory and the other a result about periodic points
for the classical $\beta$-transformation.
Finally, in
Section 6, we will make some remarks about putative elliptic
dynamical systems with the precise periodic
point behaviour and discuss possible connections with mathematical
physics.

\section{The solenoid}

Given $F(x)=bx-a$, with $a,b\in \bbZ $ coprime, let $X$ denote the 
subgroup of $\bbT ^{\Bbb Z }$ defined
by
\begin{equation}\label{sol}X=\{{\bold x}=(x_k):bx_{k+1}=ax_k\}.
\end{equation}
The group $\bbT ^{\Bbb Z }$ is compact by Tychonoff's
theorem, and $X$ is 
a closed 
subgroup so it too is compact, an example of a (1-dimensional)
{\sl solenoid}. More generally,
a solenoid is any compact, connected, abelian
group with finite topological
dimension (see
\cite{hewitt-ross}).
The automorphism $T$ is defined by the left shift-action
\begin{equation}\label{leftshift}
T({\bold x})_k=x_{k+1}.
\end{equation}
The map $T$ has the properties 1 and 2 of Section 1 by
\cite{schmidt-1995-algebraic-dynamical-systems} (see also
\cite{everest-ward-1999} for a more elementary discussion). 
In other
words, 
$$
h(T)=\log \max \{|a|,|b|\}=m(bx-a)=m(F)
$$
(a form of Abramov's formula).
Our assumption on the zero of $F$ not being
a unit root amounts to $a\neq b$, and 
the periodic points are given by
\begin{equation}\label{abperpointsformula}
|\per_{n}(T)|=|b|^n|\phi_n(a/b)|=|b^n-a^n|.
\end{equation}
At the end of this section, we will show how the periodic points 
formula (\ref{abperpointsformula})
comes about. 

In order to motivate the name, and
what follows, we will now give a second equivalent definition of
the solenoid and the action of $T$ upon it. Define $X$ to be
the topological dual of the ring $\bbZ [1/ab]$, written 
$X=\widehat {\bbZ [1/ab]}$.
Then define $T$ to be the map which is dual to the map 
$x\mapsto \frac{a}bx$ on $\bbZ [1/ab]$. The adelic point of view
arises because $X$ is isomorphic to the quotient
of $\bbR\times\prod_{p|ab}\bbQ_p$ by the diagonally embedded 
discrete subgroup $\bbZ[1/ab]$ (this is a simple finite
version of the standard adelic construction of the
dual of an $\bbA$-field, see \cite[Section 3]{chothi-everest-ward-1997}
or \cite[Chapter IV]{weil-1974-number}).
Each character on $\bbR$ restricts to a character on
$\bbZ[\frac{1}{ab}]$; this induces a map
from $\bbR\cong\hat{\bbR}$ into $X$ (injective
since $\bbZ[\frac{1}{ab}]$ is dense in $\bbR$).
The fact that the real line is `wrapped' densely into
the compact group $X$ accounts for the name solenoid.
The group $X$ is a semi-direct product of $\bbT$ by
$\prod_{p\vert ab}\bbZ_p$. The action does not
preserve the various local components, but
a direct calculation of the entropy formula is
possible (see \cite{ward-msc-thesis}).
Lind and Ward simplified
this by working with the adeles proper, which live as a covering
space to the one above. In that context, the map on each component
is simply multiplication by $a/b$. Their
approach involves tensoring the dual of $X$ with $\bbQ$ which
gives quick access
to the standard results on adeles but destroys any
periodic point behaviour (see \cite[Section 3]{lind-ward-1988}).
It is probably beneficial to keep
both points of view in mind. The elliptic system in Section 5 has the
elliptic adeles as the base group, and for the finite primes,
the local map is the local $\beta$-transformation by $a/b$. Thus
it resembles the systems defined on both the solenoid and its adelic cover.

Finally, we examine how the periodic points
formula (\ref{abperpointsformula}) comes about. This will be instructive in Section 6, when we consider
a possible elliptic analogue. Suppose $b=1$ so that $F(x)=x-a$
and consider first the
simpler case where the underlying space $X$ is the (additive) circle $\bbT$.
The map sends $x\in\bbT$ to $T_F(x)=ax \mod 1$. Thus, the points of period
$n$ are the solutions of the equation $a^nx=x$ or $(a^n-1)x=0$.
Clearly there are $|a^n-1|=|\phi_n(a)|$ solutions, which is
the division polynomial evaluated at $a$.  When the underlying space
is the solenoid $X$ as in (\ref{sol}), the map is the left shift $T$
as in (\ref{leftshift}), 
so the points of ${\bold x} \in X$ having period $n$
correspond to periodic vectors $\bold y$ of length $n$. The linear equation
generated by such a vector is of the form $C{\bold y}= \bold 0$, where
$C$ is the $n \times n$ circulant matrix on the row $(a,-1,\dots)$.
The number of solutions ${\bold y}\in\bbT^n$ of this equation, and
hence the number of periodic points, is easily verified (see
\cite[Lemma 2.3]{everest-ward-1999}) to be $|\det(C)|$. From the
well-known properties of circulants, this is equal to
$|a^n-1|=|\phi_n(a)|$.

\section{Elliptic curves}

In this section we will recall some basic results about elliptic
curves and fix the 
notation. 
A good account of elliptic curves can be found in
\cite{silverman-1986} and 
\cite{silverman-1994}; all that follows in this section can be
found in those two volumes.

Denote by $E$ an elliptic curve defined over a field $K$, and by 
$E(K)$ the group of points of $E$
having co-ordinates in $K$. When $K=\bbQ$, Mordell's theorem says 
that $E(\bbQ)$ is finitely 
generated
and the torsion-free rank is referred to simply as the
{\sl rank}.
Denote by 
$\hat h:E(\bbQ)\to\bbR $ the global canonical height 
on $E(\bbQ)$, a well known analogue of Mahler's measure 
(see (\ref{mahlermeasure})). Denote by $\lambda_p$ the local 
canonical height relative to 
the
$p$-adic
valuation.
The formula that follows gives an important decomposition of the global
height as a sum of local heights (see \cite[VIII]{silverman-1986}
and \cite[VI]{silverman-1994}):
\begin{equation}\label{globalheight}
\hat h(Q)=\sum_{p\le\infty}\lambda_p(Q),\text{ for $Q\in E(\bbQ).$}
\end{equation}
For finite $p$, whenever $Q$ has good reduction under $p$,
\begin{equation}\label{localheight}
\lambda_p(Q)=\textstyle\frac{1}{2} \log \max \{|x(Q)|_p, 1\}.
\end{equation}
In general, the local height is only defined up to the addition of a
constant. The definition (\ref{localheight}) agrees with the one in 
\cite{silverman-1988}.
In \cite{silverman-1986}, each local height is normalized by adding a constant 
to make it isomorphism-invariant.  Whether 
normalized or not, (\ref{globalheight}) still holds.

If $K=\bbR $, the 
curve
$E(\bbR)$ is isomorphic to either
$\bbT$ or $C_2\times\bbT$ 
(see \cite[V.2]{silverman-1994}).  Denote
by $E_1(\bbR)$ the connected component of the identity, which is
always isomorphic to $\bbT$. 
If $K=\bbQ_p$, the curve $E(\bbQ_p)$ can be reduced modulo $p$. 
The set of points having non-singular reduction is denoted by $E_0(\bbQ_p)$ and
the kernel of the reduction is denoted by $E_1(\bbQ_p)$. 
For odd primes $p$, there is an isomorphism 
$E_1(\bbQ _p) \stackrel{\sim}{\longrightarrow} p\bbZ p$.
This isomorphism is essentially a logarithm and it comes from the theory 
of formal groups. The situation when $p=2$ is similar; for details,
see \cite[IV]{silverman-1986}.

These isomorphisms are analogous to the one from $E_1(\bbR )$ to $\bbT $. 
The local isomorphisms for all primes $p$ play a very important role
in the development of dynamical systems because they allow  
actions on the additive local curves to be transported to the local curves proper. Consider the analogous situation in Section 1, where the immanent
group is the circle. 
When $d=1$ for example, there is an isomorphism (the logarithm) from the circle to the additive group $[0,1)$. The action on the circle really arises
from an action on $[0,1)$ which is then lifted via the logarithm to the
circle itself. In the elliptic case, the
local curve is isomorphic (via the elliptic logarithm) to an additive group. Subsequently, when we define an action 
on the $p$-adic curve, it will be one that is lifted from the additive curve.
Thus, the dynamical systems which arise when the immanent group is
the elliptic curve are exactly analogous to the case where the immanent
group is the circle (or more generally, the solenoid).

Finally, we recall the elliptic analogue of the division polynomial
$x^n-1$ on the circle. If $E$ denotes an elliptic curve defined over 
$\bbQ $ then
without loss of generality
it is defined by a generalized Weierstrass equation with integral
coefficients. There is a polynomial $\psi_n(x)$ with integer
coefficients having degree $n^2-1$ and leading coefficient $n^2$
whose zeros are precisely the $x$ co-ordinates of the points of $E$
having order dividing $n$; for details see \cite{silverman-1994}. Later,
we will consider the monic polynomial $\nu_n(x)$ of degree $n-1$, whose 
zeros are the $x$ co-ordinates of the non-identity points
in $E_1(\bbR )$ having order dividing $n$, 
\begin{equation}\label{realdivpol}
\nu_n(x)=\prod_{nQ=O \atop O\neq Q\in E_1(\bsubsubR)}(x-x(Q)).
\end{equation}
The coefficients of $\nu_n(x)$ are real algebraic numbers.

\section{The $\beta$-transformation and a $p$-adic analogue}

A comprehensive introduction to ergodic theory can be found in 
\cite{walters}. Here we just recall the definitions of 
ergodicity and entropy before examining in more detail the
$\beta$-transformation and
introducing its $p$-adic analogue.
Let $T:X\to X$ be a measure-preserving transformation on the probability
space
$(X,\mu)$.
Then $T$ is {\sl ergodic} if the only
almost-everywhere invariant sets are trivial, 
in other words
if
$\mu(T^{-1}E\Delta E)=0$ implies that $\mu(E)=0$ or 1,
where $\Delta$ is the symmetric difference.

Given two open covers $\mathcal A, 
\mathcal B$ 
of the compact topological space $X$, define
their {\sl join} to be
${\mathcal{A}} \vee {\mathcal{B}} = 
	\{ A\cap B \mid A\in{\mathcal{A}}, {B\in\mathcal{B}} \}$,
and define the entropy of $\mathcal A$ to be
$H({\mathcal A})=\log N(\mathcal A)$ where $N(\mathcal A)$ is the number of
sets in a finite subcover with minimal cardinality.
The {\sl topological entropy} of a continuous map
$T:X\to X$ is defined to be
$$
h(T)=\sup_{\mathcal A}\lim_{n\to\infty}
\frac{1}{n}H\left(\bigvee_{j=0}^{n-1}T^{-j}(\mathcal A)\right),
$$
where the supremum is taken over all open covers of $X$
(see \cite{adler-konheim-mcandrew-1965}; the topological
entropy is a measure of orbit complexity introduced as an
analogue of the measure-theoretic entropy).

The $\beta$-transformation $T_{\beta}$ is defined for real $\beta>0$ 
on the interval $[0,1)$ 
by $T_{\beta}(x)=\{\beta x \}=\beta x \pmod 1$. 
If $\beta >1$, the $\beta$-transformation
preserves an absolutely continuous probability measure
with respect to which it is ergodic \cite{renyi-1957},
the (measure--theoretic and topological)
entropy is $h(T_{\beta})=\log \beta$ (see \cite{parry-1960} and 
\cite{parry-1964}, \cite{hofbauer-1978})
and (see \cite{flatto-lagarias-1994})
the asymptotic growth rate of the periodic points 
equals
the entropy. The result about the asymptotic growth rate
will be applied in Section 5, (see (\ref{perbeta})). Strictly speaking, the definition of
topological entropy in terms of open covers does not apply
to the classical $\beta$-transformation because it has a
discontinuity; the topological entropy referred to is
that of an associated shift system (see \cite[Section 7.3]{walters}).
If $\beta \le 1$, the map is simply multiplication by $\beta$.
If $\beta <1$, $T_{\beta}$ does not preserve an
absolutely continuous measure, it has topological entropy zero 
and has no periodic points apart from 0. In all cases, the entropy
is $h(T_{\beta})=\log ^+ \beta$.

Now we define a $p$-adic analogue of the $\beta$-transformation.
For any $q\in\bbQ_p$, define a map denoted $T_q$, sometimes 
referred to as the $q$-transformation, as follows. Let $x$ be a 
generic element of 
$\bbZ_p$
and write $qx=\sum _{i=m}^{\infty} b_i p^i$. Define
$$
T_q(x)=\sum_{i=\max\{0,m\}}^{\infty}b_i p^i.
$$

In other words, $T_q$ multiplies by $q$ and cuts away the
fractional tail in order to come back to $\bbZ_p$.
Note that $T_q$ could be defined over $p\bbZ_p$ in an analogous way, 
and the ergodic
properties would not change once the Haar measure had been
normalized again.

\begin{enumerate}
\item If $|q|_p\ge1$, the map $T_q$ preserves Haar measure
on $\bbZ_p$.
\item If $|q|_p<1$ then $T_q$ is multiplication 
by $q$, and it only preserves the point mass at the
identity. 
\item The ring of $p$-adic integers $\bbZ_p$ is
homeomorphic
to the space $X= \prod_{n\in \Bbb N} Y$ of one-sided 
sequences with 
elements in $Y= \{0, \ldots, p-1 \}$, and $T_{1/p}$ is conjugate to 
the left shift $\sigma$ on $X$.
\end{enumerate}
\begin{theorem}\label{entropy}
The topological entropy of the $p$-adic
$q$-transformation is given by 
$h(T_q)=\log^+ |q|_p$.
\end{theorem}
\begin{pf}
We follow Bowen \cite{bowen-1971} and compute the
topological entropy as a volume growth rate.
It is a straightforward computation to check that that Haar measure 
on $\bbZ_p$ is
$T_q$-homogeneous, so (see \cite[Proposition 7]{bowen-1971})
\begin{equation}\label{bowen}
h(T_q)=
\lim_{m\to\infty}
\limsup_{n\to\infty}-\frac{1}{n}
\log\mu\left(
\bigcap_{k=0}^{n-1}T_q^{-k}B_m\right)
\end{equation}
where $B_m=p^m\bbZ_p$.

If $\vert q\vert_p\le 1$,
$T_q^{-1}B_m\supset B_m$ so
$$
\bigcap_{k=0}^{n-1}T_q^{-k}B_m=B_m,
$$
and (\ref{bowen}) gives $h(T_q)=0=\log^{+}\vert q\vert_p.$

If $\vert q\vert_p=p^r>1$, then
$T_q^{-1}B_m=B_{m+r}$, so
$$
\bigcap_{k=0}^{n-1}T_q^{-k}B_m=B_{m+rn},
$$
so by (\ref{bowen}) $h(T_q)=r\log p=\log^{+}\vert q\vert_p$.
\end{pf}

\begin{theorem}\label{ergodicity}
Let $q\in\bbQ _p$ with $|q|_p\ge 1$. 
The map $T_q$ is ergodic with respect to
Haar measure for $|q|_p>1$, and is
not ergodic for $\vert q\vert_p=1$.
\end{theorem}

\begin{pf}
Assume $|q|_p>1$ and let $\mathcal A$ denote the algebra of all 
finite unions of measurable
rectangles and suppose $E$ is a measurable set invariant under $T_q$.
For any given $\epsilon>0$ it is possible to choose $A\in \mathcal A$ with 
$\mu (E \bigtriangleup A)<\epsilon$, and thus
$|\mu(E)-\mu(A)|<\epsilon.$
Choose $n$ such that $B=T_q^{-n} A$
depends upon different co-ordinates  from $A$: then 
$\mu(A\cap B)=\mu(A)\mu(B)=(\mu(A))^2$.
Also, $\mu(E \bigtriangleup B)<\epsilon$, and 
$\mu(E\bigtriangleup
 (A\cap B)) \le \mu((E\bigtriangleup
 A)\cup (E\bigtriangleup
B))<2\epsilon$.
Thus, $|\mu(E) - \mu(A\cap B)|<2\epsilon$ and 
$|\mu(E)-\mu(E)^2| < 4\epsilon. $			
Since $\epsilon$ is arbitrary, this implies $\mu(E)=0$ or $1$ and thus $T$ 
is
ergodic.

When $q$ is a unit, the open sets of the form
$p^n\bbZ_p$ for $n\ge 1$ are all
invariant under $T_q$.
\end{pf}
\begin{remark} Notice that in this setting ergodicity and mixing coincide. 
Coelho and Parry have studied the ergodic decomposition
of $T_q$ when $q$ is a unit (see \cite{coelho-parry-1998}).
\end{remark}

A consequence of properties 1 and 2 for the systems in
Section 1 is that the logarithmic growth rate of the
periodic points coincides with the entropy (see
\cite{everest-ward-1999}). That this also holds for
$T_q$ follows from the next result.

\begin{theorem}\label{periodicpoints} Given $q\in\bbQ_p
\setminus U$, where $U$ denotes the
set of unit roots in $\bbQ_p$, let 
$T_q$
denote the $q$-transformation on $\bbZ_p$.
Then \begin{equation}\label{countper}
\log |\per_{n}(T_q)|=n\log ^+|q|_p.
\end{equation}
\end{theorem}

\begin{pf}
Firstly, consider the case  $|q|_p<1$.
Then for $n\to\infty$
$T_q^n(x) \rightarrow 0$ for all $x\in\bbZ_p$.
Thus 
$T_q$
has only one periodic point (zero) and both sides of (\ref{countper})
are zero.

When  $|q|=1$,
the action of $q$ on $\bbZ_p$ is simply multiplication, so
the periodic points are solutions to the equation
$q^nx=x$. Since $q$ is not a unit root, there are
no periodic points except $x=0$, so
(\ref{countper}) holds.

Finally suppose $|q|_p>1$. If $q=p^{-k}$ with $k>0$, the periodic 
points are easy to determine. 
We  have 
$T_q^n(x)=\sum_{i=0}^{\infty} a_{i+nk} p^i$ and the solutions to  
$T_q^nx=x$ are
given by the $p^{kn}$ points with $a_{i+nk}=a_i$ for $i=0, \dots, kn-1$.
Thus, both sides of (\ref{countper}) are equal to $nk\log p$. 
In general, suppose  $|q|_p=p^k$. We claim that for each integer
$a$ with $0\le a <p^{nk}$, there is a unique $y\in \bbZ _p$ with
$T_q^n(a+p^{nk}y)=a+p^{nk}y$. This follows because the left hand side
is $b+q^np^{nk}y$ for some $b\in \bbZ_p$, which depends only upon $a,q$
and $n$. Write $q^np^{nk}=v$ for some $p$-adic unit $v$ then the equation 
$b+vy=a+p^{nk}y$ has a unique solution for $y\in \bbZ_p$.
This shows that there are at least $p^{nk}$ solutions of $T_q^nx=x$.
That there can be no more follows because we may take the $a$ as
above as coset representatives for $\bbZ_p/p^{nk}\bbZ_p$ so
every element $x\in \bbZ_p$ is represented by some $a$.
\end{pf}

In conclusion, given any $p$-adic elliptic curve $E$ and any
point $Q \in E(\bbQ_p)$, we can construct a dynamical system
in the following way. The curve is locally isomorphic to the
group $p\bbZ_p$ and therefore to $\bbZ_p$. Now let $Q$ act 
via the $q$-transformation  
on the additive curve, where $q=x(Q)$. Then transport this action to the curve proper
via the logarithm. This is an exact analogue of the
toral dynamical systems in Sections 1 and 2.

\section{Dynamics on the Elliptic Adeles}

From here on, 
let $E$ denote an elliptic curve defined over $\bbQ$ and let
$Q\in E(\bbQ)$. 
The explicit formula (\ref{localheight}) for the local height of $Q$
does not hold if $Q$ has bad
reduction or if $p$ is the prime at infinity. In
particular, the local height in these cases can be negative. Since
the entropy of a map is never negative, we will work with points whose
local heights are guaranteed to be non-negative.

\begin{claim}
There exists an $n\ge1$
for which the finite-index subgroup
$nE(\bbQ )\le E(\bbQ )$
has $\lambda_p(Q) \ge 0$ for all $p<\infty$ and 
$Q\in nE(\bbQ )$.
\end{claim}

\begin{pf} 
Since $E_1(\bbQ )$ is a subgroup of finite index in $E(\bbQ)$ 
(\cite[VII]{silverman-1986}), for each bad prime $p$ there exists an integer
$n_p$ such that $E_1(\bbQ )$ has index $n_p$ in $E(\bbQ)$. 
Let $n=\prod_{\text {bad }p} n_p$, then $nE(\bbQ)\le
E_1(\bbQ)$ 
for all bad $p$.
Recall now that if $Q$ has good reduction at $p$ 
(and this includes the case where $Q\in E_1(\bbQ )$) then the local
height at $p$ is given by (\ref{localheight}) and it follows that 
$\lambda_p(Q) \ge 0$.
\end{pf}

Define $\cS$ to be the set of bad primes together with infinity. Assume 
that the point
$Q$ satisfies
\begin{equation}\label{assumptions}
\lambda_p(Q) > 0 \quad \textrm{ for all } p\in \cS.
\end{equation}
If $Q\in nE(\bbQ)$ then $Q\in E_1(\bbQ)$ for
all the bad primes. It follows from (\ref{localheight}) that the
local height is actually positive. If the rank
of $E(\bbQ)$ is not zero then $nE(\bbQ)$ has finite index in
$E(\bbQ)$ so
in that case, there is a large stock of points
$Q$ which satisfy (\ref{assumptions}). At the infinite prime, this
amounts to
assuming
that $Q$ lies in a neighbourhood of the identity.

Suppose $Q\in 
E(\bbQ )$ 
is a point 
for which the assumption (\ref{assumptions}) holds. Define $X$ 
to be the space
\begin{equation}\label{defineX}
X=\prod_{p\le\infty}E_1(\bbQ_p).
\end{equation}
The point $Q$ induces an action $T_Q:X\rightarrow X$ in the
following way: $(T_Q)_p$ is 
the $q$-transformation if $p$ is finite (where $q=x(Q)$)
and the $\beta$-transformation if $p$ is infinite, where 
$\log \beta=2\lambda_{\infty}(Q)$. Remember that these are actions on
$\bbT $ and $p\bbZ _p$, but for every $p$, the action can be
transported to $E_1(\bbQ _p)$ via the isomorphisms in Section 3.
The statements in the following theorem are analogues of
statements 1 and 2 in the introduction. There we supposed
that the zeros of $F$ were not torsion points of $\bbT $. The assumption
that $Q$ is not a torsion point of $E$ is built into
(\ref{assumptions}):
$Q$ is a torsion point if and only if $\hat h(Q)=0$ and
(\ref{assumptions})
guarantees that $\hat h(Q)>0$.

\begin{theorem}\label{teo1}
With the definitions and assumptions above,
\begin{enumerate}
\item the entropy of $T_Q$ is given by $h(T_Q)=2\hat h(Q)$ and 
\item the asymptotic growth rate of the periodic points is given  
by the division polynomial $\nu_n(x)$ in
{\rm(\ref{realdivpol})}: 
$$\log|\per_{n}(T_Q)|\sim\log |b^n\nu_n(q)|\text{ as } 
n\rightarrow\infty.$$ 
\end{enumerate}\end{theorem}

\begin{pf}
By Theorem \ref{entropy}, the entropy of each component of $T_Q$ is 
given by
$\log\beta_p$, where $\beta_p=\beta$ if $p=\infty$ and 
$\beta_p=\max \{|x(Q)|_p,1\}$ if $p$ is finite. Since there are only 
finitely many primes for which the
local dynamical systems are not isometries, 
Theorem 4.23 in \cite{walters} applies giving 
$$
h(T_Q)=h(T_{\beta})+\sum_{p<\infty} h(T_q)=
\sum_{p\le \infty} \log\beta_p=
2\sum_p\lambda_p(Q)=2\hat h(Q). $$
For the asymptotic growth rate of the periodic points note that
if dynamical systems 
$\hat T_i:X_i\to X_i$ ($i=1,\dots,r$)
are given and the point $x_i$ has period $m$ under $\hat T_i$ for 
$i=1,\dots ,r$ 
then $(x_i)$
has period $m$ under $\prod_i\hat T_i$. Thus we may count the 
contribution
to the periodic points from each prime separately. For $p<\infty$, 
from Theorem \ref{periodicpoints}(\ref{countper})
\begin{equation}\label{finiteperiodic}
\log |\per_{n}(T_q)|=n\log ^+|q|_p=-n\log |b|_p.
\end{equation}
Note that our assumption on $Q$ guarantees that $q$ is not
an integer and so, in particular, $q$ is not a root of unity.
Summing over all finite $p$ and using the product formula, we
obtain a total contribution of $n\log |b|$ to the periodic points.
For the infinite prime, we quote a deep result from
\cite{flatto-lagarias-1994} 
which says
\begin{equation}\label{perbeta}
\log |\per_{n}(T_{\beta})|=n\log \beta+o(n).
\end{equation}
From (\ref{finiteperiodic}) and (\ref{perbeta}) we have the 
formula
\begin{equation}\label{perpointsform}
\log|\per_{n}(T_Q)|=n\log |b| + n\log \beta + o(n).
\end{equation}
Finally, we quote from Theorem 6.24 in \cite{everest-ward-1999} which gives
\begin{equation}\label{perelliptic}
\log |\nu_n(q)|=n\log \beta + o(n).
\end{equation}
The formula (\ref{perelliptic}) depends upon an application
of the elliptic analogue of Baker's theorem from
transcendence theory (see \cite{david-1991},
\cite[Section 7]{everest-ward-1998}
and \cite[Theorem 6.18]{everest-ward-1999}).
It follows from (\ref{perpointsform}) and
(\ref{perelliptic}) that $\log |b^n\nu_n(q)|$
is asymptotically equivalent to $\log |\per_{n}(T_Q)|$.
\end{pf}

Set now $\cS^*(Q)=\{p : |x(Q)|_p > 1\} \cup \{\infty\}$, let 
$X^*=\prod_{p\in \cS^*(Q)}$ and let $T_Q^*$ be defined component-wise as 
above.
\begin{theorem}\label{teo2}
\begin{enumerate}
\item The entropy of $T_Q^*$ is given by $h(T_Q^*)=2\hat h(Q)$, 
\item the asymptotic growth rate of the periodic points is given  
by the division polynomial {\rm(\ref{realdivpol})}: 
$$\log|\per_{n}(T_Q^*)| \sim \log |b^n\nu_n(q)|,$$ 
\item $T_Q^*$ is ergodic.
\end{enumerate} \end{theorem}

\begin{pf}
For the entropy and the periodic points, the same arguments as in 
Theorem \ref{teo1} holds giving the
desired result. The ergodicity is proved in Theorem \ref{ergodicity}
\end{pf}

The pros and cons of our construction may be summarized as follows.
Firstly, we have constructed a dynamical system whose immanent group
is the adelic elliptic curve. The map is defined locally by the
$p$-adic $\beta$-transformation on the additive curve.
Secondly, the construction exhibits phenomena which resemble those
in the solenoid case. Against these comments we must set the
following. Firstly, the maps we are using are not continuous because of the discontinuity of the classical 
$\beta$-transformation. The effect upon the map $T_Q$ is to deny
continuity at infinity on the archimedean component. Secondly, we would have preferred to see periodic 
point
behaviour which was counted
precisely by the usual elliptic division
polynomial (rather than just asymptotically by the real division polynomial). Thirdly, the map at the
archimedean prime uses {\it a priori} knowledge of the archimedean height
of the point. Fourthly, we made special assumptions to guarantee that
each local height was non-negative. Although these assumptions were natural, 
at the
infinite prime and each
bad prime we assumed our point was to be found in a neighbourhood of
the identity, we would have preferred not to have needed any assumptions. 
In the next section, we will discuss how these
deficiencies might be overcome.

\section{Putative Elliptic Dynamics}

Suppose $E$ denotes an elliptic curve defined by a generalized
Weierstrass equation with integral coefficients. For each $n\in\bbN$,
let $\psi_n(x)=n^2x^{n^2-1}+...$ denote the $n$-th division polynomial.
Let $Q$ denote a non-torsion rational point on $E$, with $x(Q)=a/b$.
It is tempting to
conjecture that
there must a compact space $X$ with a continuous action 
$T_Q:X\rightarrow X$ whose entropy of is given by 
$h(T_Q)=2\hat h(Q)$, and whose periodic points are counted, 
in the sense of Section 1, by $E_n(Q)=|b^{n^2-1}\psi_n(a/b)|$. 
The sequence $E_n(Q)$ is
certainly a divisibility sequence like its toral counterpart. 

However, there is a problem in making
the obvious conjecture. Hitherto, the systems we have considered have 
been $\bbZ$-actions
(see \cite{everest-ward-1999} or \cite{schmidt-1995-algebraic-dynamical-systems} for the definition of $\bbZ^d$-action).
Recent work (see \cite{puri-ward-1998}) makes it unlikely that
the sequence $E_n(Q)$  counts periodic points for
a $\bbZ$-action. Indeed, when $E$ is given by the equation
$y^2+y=x^3-x$ and $Q$ is the point $Q=(0,0)$, it follows from
\cite{puri-ward-1998} that $E_n(Q)$ cannot represent periodic point data for any $\bbZ$-action. The sequence begins 1,1,1,1,2\dots
so it violates the divisibility condition (with $n=5$) given in
(3) of \cite{puri-ward-1998}. 
What does seem possible is that $E_n(Q)$ represents periodic
point data for a $\bbZ^2$-action. The two reasons for saying this
are firstly that the growth rate of the sequence is quadratic exponential
in $n$ (see \cite[Theorem 6.18]{everest-ward-1999}). Thus, the
sequence is more likely to represent
$|\per_{(n\bsubZ)^2}(T)|$, for some $\bbZ^2$-action $T$, where
$(n\bbZ)^2$ represents the subgroup $\bbZ^2$ having index $n^2$ 
and consisting
of all $(x,y)$ with $n|x$ and $n|y$. This would make it consistent with
the known properties of algebraic $\bbZ^2$-actions 
(see \cite[Theorem 7.1]{lind-schmidt-ward-1990}). Secondly, 
the general
feeling persists that natural maps on elliptic curves tend to
be quadratic.  

We conjecture that for every rational elliptic curve $E$
and every rational point $Q \in E(\bbQ)$, there 
is a (necessarily infinite dimensional) compact space $X$ with a 
continuous $\bbZ^2$-action 
$T_Q:X\rightarrow X$ having the
following properties:
\begin{enumerate}
\item the entropy of $T_Q$ is given by $h(T_Q)=2\hat h(Q)$, and
\item the periodic points are counted by $\psi_n$ in the sense that
$$
|\per_{(n\Bbb Z)^2}(T_Q)|=|b^{n^2-1}\psi_n(a/b)|.
$$
\end{enumerate}

The integral case (where $b=1$) of the conjecture is already challenging. Suppose 
then $b=1$ and we seek
an action of the integral point $Q$, where $x(Q)=a$. We would hope to recognize the sequence $|\psi_n(a)|$
in some natural way as counting periodic points. 
In \cite[VI.4]{everest-ward-1999},
we noted that the numbers  $|\psi_n(a)|$ arise as determinants of 
a nested sequence of integral $n \times n$ Hankel
matrices. We suggest that these matrices might be the analogues of
the circulant matrices $C$ in Section 2.

We hope that an action with such beautiful dynamical data would not
be deficient in the way that the map of Section 5 was deficient.
In particular, the potential negativity of the local heights need
not seem such a threat. Although a negative entropy cannot exist,
nonetheless, the difference between two non-negative entropies
can make sense. If one dynamical systems extends another then the
difference between their two entropies represents the entropy across
the fibres. This raises the possibility that a phenomenon such as
bad reduction might well have a dynamical interpretation.

Interest in our conjecture (see \cite[Question 14]{everest-ward-1999})
is heightened because of the connection with the following remarkable
circle of ideas. On the one hand, mathematical
physicists have studied
the dynamics of integrable systems
(see \cite{veselov-1991},
\cite{veselov-1992}). Here, {\it inter alia},
one looks for meromorphic maps on the complex plane which
commute with polynomials. It is a classical result of
Ritt (see \cite{ritt-1923}) that all non-trivial examples arise from
the exponential function or the elliptic functions associated 
to some lattice. Coincidentally, Morgan Ward
(see \cite{mward-1948b}, \cite{mward-1948a})
showed that all integer sequences satisfying a certain
natural recurrence
relation arise from the exponential function or the elliptic
functions associated to some lattice, suitably evaluated.
In the exponential case, these sequences can always be identified
with the periodic point data for toral automorphisms. It
is hoped
that in the elliptic case also, the sequences $|\psi_n(a)|$
represent the periodic point data for some elliptic systems.
That being so, a new chapter in integrable systems could be
written, yielding further inter-play between elliptic curves
and mathematical physics.

\bibliographystyle{plain}

\end{document}